\documentclass[10pt,a4paper]{article}

\usepackage[margin=0.85in]{geometry}
\usepackage{graphicx}
\usepackage{amsmath,amssymb}
\usepackage{booktabs}
\usepackage{bm}
\usepackage{cite}
\usepackage{hyperref}
\hypersetup{hidelinks}
\usepackage{microtype}
\newcommand{\authorchinesename}{\raisebox{-0.15ex}{\includegraphics[height=1.2em]{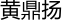}}}
\newcommand{\pubcite}[1]{\textup{[#1]}}

\title{Geometric priors for single-patch neural Calabi--Yau metrics on the quintic: a conservative fixed-$\psi$ Ablation study}
\author{D~Yang Eng (\authorchinesename)\textsuperscript{1,2,*}\\[0.5em]
\small \textsuperscript{1}Technical University of Munich (TUM), Arcisstraße 21, 80333 Munich, Germany\\
\small \textsuperscript{2}Singapore Institute of Technology, 1 Punggol Coast Road, Singapore 828608, Singapore\\[0.4em]
\small \textsuperscript{*}Author to whom any correspondence should be addressed.\\
\small\texttt{dyang.eng@tum.de}
}
\date{}

\begin{document}
\maketitle
\begin{abstract}
We investigate how geometric priors shape neural approximations of the complex Monge--Amp\`ere equation on the quintic Calabi--Yau threefold. This fixed-$\psi$ diagnostic compares four single-patch ans\"atze: a coordinate SIREN potential model, a direct-metric Cholesky baseline, an algebraic section ansatz, and a section-plus-neural hybrid. In a 96-run learning-rate and frequency archive, section-based models reach validation $\sigma=4.14$--$4.16\%$, whereas SIREN and Cholesky remain at $41.29\%$ and $44.27\%$. Multi-seed controls use a shared affine-$\eta$ evaluator and five fresh validation clouds. Stable section settings give $3.84$--$4.77\%$, compared with $28.46$--$30.59\%$ for Cholesky and $34.73$--$38.13\%$ for SIREN. Hybrid projective features do not improve best-case accuracy, but they stabilise the stressed regime. They limit $\sigma$ to $4.26$--$4.41\%$, versus $4.45$--$205.45\%$ for raw coordinates. Within this protocol, section-based geometric structure is the main driver of accuracy; coordinate-only models plateau at substantially higher error, and raw hybrid flexibility can amplify instability. We provide a projective-scaling derivation, document protocol asymmetries, and identify projectively consistent hybrid corrections as the next step.
\end{abstract}

\par\bigskip

\section{Introduction}
\label{sec:introduction}

The problem of constructing Calabi--Yau metrics originates in Calabi's conjecture. Calabi proposed that, on a compact K\"ahler manifold, a representative of a fixed K\"ahler class should exist with prescribed Ricci form; in the case $c_1(X)=0$, this predicts a unique Ricci-flat K\"ahler metric in each K\"ahler class~\pubcite{1, 2}. Yau proved the conjecture by solving a complex Monge--Amp\`ere equation~\pubcite{3, 4}. The theorem is fundamental for both differential geometry and string compactification, but it is nonconstructive: for compact Calabi--Yau threefolds such as the quintic, explicit closed-form Ricci-flat metrics are not known.

This gap between existence and explicit computation has persisted for nearly five decades. It motivates both classical balanced-metric algorithms and, more recently, neural approximations to Calabi--Yau metrics. Physics-informed neural networks (PINNs), neural operators, and related scientific machine-learning models have become practical tools for approximating nonlinear partial differential equations and expensive numerical solvers~\pubcite{5--7}. Their success, however, is rarely a consequence of raw network capacity alone. Across recent scientific-ML work, performance depends strongly on how the architecture represents the underlying geometry and how physical constraints enter the loss~\pubcite{8--11}. These issues become sharper for geometric PDEs, where the learned object is not a scalar field but a metric, and downstream quantities depend on its determinant, inverse, and derivatives.

This paper studies these issues in a controlled setting: the complex Monge--Amp\`ere equation for Ricci-flat K\"ahler metrics on the quintic Calabi--Yau threefold. The aim is diagnostic rather than competitive. We do not claim a new state-of-the-art Calabi--Yau metric approximation, a complete phenomenological pipeline, or a universal limitation of neural solvers. Instead, we ask a narrower question: when persistent plateaus appear in a fixed single-patch protocol, are they better explained by hyperparameter tuning, raw model capacity, geometric constraint enforcement, or the representation built into the ansatz?

\subsection{Background for Machine-Learning Readers}
\label{sec:ml-background}

The geometric terminology used in this paper can be translated into the language of constrained function approximation. A Calabi--Yau threefold is a compact complex three-dimensional K\"ahler manifold with vanishing first Chern class, $c_1(X)=0$. The quintic hypersurface is a standard example. After fixing a K\"ahler class, the Calabi--Yau theorem guarantees a unique Ricci-flat K\"ahler metric in that class. Here ``Ricci-flat'' means that the Ricci curvature tensor of the metric vanishes. In string compactification, Calabi--Yau backgrounds solve the leading supersymmetric vacuum conditions~\cite{candelas1985}.

For a K\"ahler metric, the metric tensor is generated locally from a scalar potential $K$ by $g_{i\bar{j}}=\partial_i\partial_{\bar{j}}K$. The complex Monge--Amp\`ere equation is the nonlinear PDE that enforces Ricci-flatness at the level of the volume form: the metric volume $\omega^3$ must be proportional to the canonical holomorphic volume form $\Omega\wedge\bar{\Omega}$. Numerically, one evaluates a local density ratio and checks whether it is constant over sampled points. The statistic $\sigma$, also denoted $\sigma_{\mathrm{RMS}}$, is the weighted root-mean-square relative deviation of this ratio from its mean; smaller values indicate a closer Monge--Amp\`ere solution. The Hodge numbers $h^{1,1}$ and $h^{2,1}$ count, respectively, K\"ahler and complex-structure deformation directions. For the quintic, $h^{1,1}=1$ and $h^{2,1}=101$. Table~\ref{tab:notation} summarises the notation used in the rest of the manuscript.

\begin{table}[t]
\centering
\caption{Main notation used in the manuscript. The table is included to make the geometric and learning notation explicit for readers from machine learning and numerical PDEs.}
\label{tab:notation}
\begin{tabular}{p{0.24\linewidth}p{0.68\linewidth}}
\hline
Symbol & Meaning \\
\hline
$X_{\psi}$ & Quintic Calabi--Yau hypersurface in $\mathbb{CP}^4$ at complex-structure parameter $\psi$. \\
$z$, $w$ & Homogeneous projective coordinates $z$ and affine chart coordinates $w_i=z_i/z_0$. \\
$K$ & Total K\"ahler potential used to generate a metric. \\
$f_{\theta}$ & Neural scalar correction or coordinate-network output, with trainable parameters $\theta$. \\
$\omega_0$ & SIREN frequency scale swept in the ablation grid. \\
$g_{i\bar{j}}$ & Hermitian metric tensor, usually computed as $\partial_i\partial_{\bar{j}}K$. \\
$\Omega$ & Holomorphic $(3,0)$-form defining the Calabi--Yau volume reference. \\
$\eta$, $\rho$ & Local metric-density ratios used for training and validation. \\
$\kappa$ & Weighted normalisation constant for the density ratio. \\
$\sigma$ or $\sigma_{\mathrm{RMS}}$ & Weighted RMS relative dispersion of the normalised density; reported as a percentage. \\
$s(z)$, $H$ & Degree-three section vector and learned Hermitian matrix in $K_{\mathrm{sec}}=\log(s^\dagger Hs)$. \\
$N_{\mathrm{val}}$ & Number of validation points used to estimate $\sigma$. \\
\hline
\end{tabular}
\end{table}

\subsection{Geometric Setting}

The formal geometric setting is as follows. Yau's theorem guarantees that every K\"ahler class on a compact Calabi--Yau manifold contains a unique Ricci-flat K\"ahler metric~\pubcite{3, 4}. In terms of a K\"ahler potential $K$, this metric is governed by the complex Monge--Amp\`ere equation
\begin{equation}
    \left(\frac{i}{2}\partial\bar{\partial}K\right)^3
    = c\,\Omega\wedge\bar{\Omega},
\end{equation}
where $\Omega$ is the holomorphic $(3,0)$-form and $c>0$ is fixed by the total volume. The theorem is existential rather than constructive, so explicit compact examples must be approximated numerically.

The classical finite-dimensional approach represents $K$ through Hermitian forms on spaces of holomorphic sections. Donaldson's balanced-metric programme~\pubcite{13--15}, its implementation on the quintic by Douglas et al~\cite{douglas2008}, the K3 numerical solver of Headrick and Wiseman~\cite{headrick2005}, and the energy-functional refinements of Headrick and Nassar~\cite{headrick2013} provide the main algebraic reference points. Bergman-kernel asymptotics motivate why increasing the section degree improves the approximation~\pubcite{19, 20}; the rapid growth of the finite section space is then a direct combinatorial consequence of the degree-$k$ homogeneous monomial basis.

Neural approaches offer a complementary route. Early work trained networks on Donaldson-generated data~\cite{ashmore2020}; holomorphic-network methods showed the value of section-aware architectures~\cite{douglas2022}; and later work trained directly from Monge--Amp\`ere or related residuals, incorporated moduli dependence, and developed packages such as \texttt{cymetric}, \texttt{CYJAX}, and \texttt{cymyc}~\pubcite{23--27}. Recent symmetry-aware, symbolic, and section-based methods further show that explicit geometric structure can improve accuracy, stability, or interpretability~\pubcite{28--30}. In particular, Grassmannian learning keeps Donaldson's section framework but uses optimisation to select efficient section subspaces and Hermitian data~\cite{ek2024}. The present paper is positioned inside this line of work, but it isolates a smaller question than those full pipelines address.

\subsection{Benchmark and Gap}

We work with the standard Dwork one-parameter family of quintic hypersurfaces in $\mathbb{CP}^4$~\cite{candelas1991},
\begin{equation}
    X_{\psi}
    =
    \left\{
        z \in \mathbb{CP}^4 :
        \sum_{i=0}^{4} z_i^5 - 5\psi \prod_{i=0}^{4} z_i = 0
    \right\},
\end{equation}
and fix the complex-structure modulus to the Fermat point $\psi=0$. The quintic is a useful diagnostic benchmark because its topology is well understood, with $h^{1,1}=1$ and $h^{2,1}=101$~\pubcite{12, 31}; because algebraic and spectral numerical baselines are available~\pubcite{15, 16, 18, 32}; because its Monge--Amp\`ere formulation imposes hard geometric constraints; and because it appears repeatedly in neural Calabi--Yau studies~\pubcite{23--25, 27}.

The difficulty is also representative. The Monge--Amp\`ere residual depends on second derivatives of $K$, while curvature-sensitive quantities, Yukawa couplings, and Kaluza--Klein spectra require still stronger control of derivatives, inverse metrics, and global consistency~\pubcite{27, 33--35}. Published neural studies use different losses, sampling procedures, batch sizes, geometric constraints, and random seeds. As a result, cross-paper numbers alone do not identify whether an improvement is caused by architecture, optimisation, data generation, or evaluation convention. This paper addresses that ambiguity by comparing four architectures inside one archived workflow and then checking the ordering under a smaller shared-evaluator control suite.

\subsection{This Work}

The main archive is a full-factorial $4\times4\times6=96$-run fixed-$\psi$ ablation. It spans four model architectures, four learning rates $\{10^{-5},10^{-4},10^{-3},10^{-2}\}$, and six SIREN frequency parameters $\omega_0$. The architectures are chosen to form a ladder of geometric structure:
\begin{itemize}
    \item \textbf{SIREN potential}: a coordinate network outputs a scalar K\"ahler potential.
    \item \textbf{Cholesky metric}: a coordinate network predicts a positive-definite Hermitian metric directly.
    \item \textbf{Section-Pure}: the K\"ahler potential is $K_{\mathrm{sec}}=\log(s(z)^\dagger Hs(z))$, where $s(z)$ contains all $35$ degree-three monomial sections of $\mathcal{O}(3)|_X$ and $H$ is learned.
    \item \textbf{Section-Hybrid}: a small SIREN correction is added to the section potential, $K_{\mathrm{hyb}}=K_{\mathrm{sec}}+0.1f_\theta(z)$.
\end{itemize}

This design separates three effects that are often entangled: coordinate learning, positivity by direct metric parametrisation, algebraic section structure, and extra neural flexibility on top of that structure. The archive is supplemented by a multi-seed control suite using seeds $\{11,22,33\}$, matched batch sizes $\{2000,4096\}$, and a shared post-hoc affine-$\eta$ evaluator. The control suite is not a full benchmark replacement; it tests whether the architecture ordering seen in the archive survives the most obvious seed, batch-size, and reporting asymmetries.

The result is simple but informative. Section-based models reach the $4\%$--$5\%$ $\sigma$ regime, while the coordinate-only baselines remain far higher: roughly $30\%$--$44\%$ depending on whether one reads the archived or controlled comparison. The identity degree-three section control already reaches $\sigma=4.44\%$, close to the trained section-Pure and low-learning-rate hybrid models. Thus, in this protocol, the dominant gain comes from the section prior itself, not from additional raw neural flexibility. The hybrid branch gives the nominally lowest archived run, but the gap to section-Pure is too small to rank the two section models statistically; its clearer effect is to introduce a stability hazard unless its features are made projectively consistent.

\subsection{Contributions}

The paper makes three contributions.
\begin{enumerate}
    \item It gives a self-contained formulation of the fixed-$\psi$ quintic Monge--Amp\`ere problem, including the affine density, induced metric, and $\sigma$ diagnostic used in the experiments.
    \item It derives the projective consistency of the degree-three section ansatz and clarifies why the raw hybrid correction is only a gauge-fixed empirical extension of that ansatz.
    \item It reports a conservative architectural diagnosis based on a complete 96-run archive and targeted multi-seed controls, using regenerated Donaldson references produced with \texttt{CYJAX}~\cite{gerdes2023cyjax}.
\end{enumerate}

The scope is deliberately narrow. The paper does not rank all Calabi--Yau metric solvers, does not replace atlas-based or spectral benchmarks, and does not claim phenomenology-grade precision. It shows that, in the fixed single-patch regime studied here, geometric section structure is the principal variable separating the observed outcomes.

\subsection{Paper Organisation}

Section~\ref{sec:related} situates the study within numerical Calabi--Yau metrics and scientific machine learning. Section~\ref{sec:methodology} defines the metric density, model families, losses, and protocol. Section~\ref{sec:results} reports the archive and control experiments. Section~\ref{sec:discussion} interprets the evidence and its limitations. Section~\ref{sec:conclusion} summarises the conclusions and next steps.

\section{Background and Related Work}
\label{sec:related}

\subsection{Finite-Dimensional Calabi--Yau Metrics}

The numerical approximation of compact Ricci-flat Calabi--Yau metrics is usually organised around finite-dimensional geometric ans\"atze. Bergman-kernel asymptotics explain why spaces of holomorphic sections recover the K\"ahler geometry in the large-degree limit~\pubcite{19, 20}. Donaldson translated this structure into a constructive approximation scheme based on balanced metrics and the T-operator~\pubcite{13--15}. In this view, increasing the section degree $k$ improves the approximation, with convergence controlled by the large-$k$ asymptotics of projective embeddings. Douglas, Karp, Lukic and Reinbacher implemented the method on the quintic and introduced the $\sigma$ diagnostic used throughout this literature~\cite{douglas2008}; Headrick and Wiseman developed an earlier K3 numerical Ricci-flat metric solver exploiting K\"ahler geometry~\cite{headrick2005}. Headrick and Nassar then showed that energy-functional optimisation within the same algebraic ansatz can produce substantially more accurate fixed-degree metrics than balanced iteration alone~\cite{headrick2013}. Spectral computations of Laplacian eigenvalues and eigenfunctions provide another classical numerical reference point for the quintic~\cite{braun2008}.

This classical background matters for the present paper in two ways. First, it explains why the section ansatz is not merely another neural feature map: it is the finite-degree truncation of a geometrically convergent family. Second, it makes comparisons delicate. A degree-$k$ balanced metric, an energy-functional optimum, a neural surrogate trained on algebraic data, and a residual-trained neural metric are not interchangeable baselines, even when they are evaluated with a similar $\sigma$ statistic.

\subsection{Neural Calabi--Yau Metrics}

\subsubsection{Progress and flexibility.}
Machine-learning approaches to Calabi--Yau metrics have developed along several lines. Ashmore, He and Ovrut demonstrated that networks can accelerate Donaldson-style metric approximation~\cite{ashmore2020}. Holomorphic-network methods introduced an explicitly section-aware neural construction for numerical Calabi--Yau metrics~\cite{douglas2022}. Anderson et al.\ introduced moduli-dependent learning for Calabi--Yau and SU(3)-structure metrics~\cite{anderson2021}. Jejjala et al trained neural metrics directly from a volume-form loss and verified discrete-symmetry invariance of the trained metrics~\cite{Jejjala2022}. The \texttt{cymetric} work extended the approach to complete intersections and Kreuzer--Skarke examples with custom losses and sampling tools~\cite{larfors2022}, while \texttt{CYJAX} provided a JAX implementation centred on algebraic K\"ahler-potential ans\"atze with automatic K\"ahlerity and patch compatibility~\cite{gerdes2023cyjax}. The \texttt{cymyc} pipeline and related curvature work connect learned metrics to Yukawas, curvature diagnostics, and precision phenomenology~\pubcite{27, 34, 36}.

\subsubsection{Structural obstructions.}
Despite these advances, the geometry of the problem imposes constraints that black-box neural architectures do not automatically respect. If a network learns the metric tensor directly, pointwise positive definiteness can be imposed by a Cholesky or related factorisation, but this does not guarantee that the resulting tensor is K\"ahler, projectively well defined, or derived from a global potential. If a network learns a potential correction, the ansatz locally produces a closed $(1,1)$ form after applying $i\partial\bar{\partial}$, but the condition that the total form, for example $\omega_{\mathrm{FS}}+i\partial\bar{\partial}\phi$, remain positive is an open pointwise constraint on a function space. Ordinary gradient steps in parameter space do not naturally preserve that positivity condition. Thus the common shorthand that a potential ansatz is ``automatically K\"ahler'' should be read carefully: closedness is automatic locally, but positivity, projective consistency, and global compatibility are additional requirements. These are not merely numerical inconveniences. They reflect a mismatch between soft, weighted-loss formulations and the hard algebraic and differential requirements of K\"ahler geometry.

The most relevant lesson from the recent literature is therefore not that one network family is universally best. Rather, the strongest pipelines encode more of the geometry: K\"ahler potentials rather than unconstrained tensors, projective or patch-compatible coordinates, moduli dependence when needed, symmetry handling, or spectral/section structure. Hendi, Larfors and Walden show that fundamental-domain canonicalisation improves accuracy on symmetric Calabi--Yau examples at fixed training budget~\cite{hendi_learning_2024}. Mirjani\'c and Mishra show that extrinsic symmetries can yield compact symbolic metric approximations on special loci~\cite{mirjanic_mishra_2025}. Ek et al take a particularly relevant route: they keep Donaldson's section framework but formulate part of the learning problem as optimisation over Grassmannian section subspaces and Hermitian data~\cite{ek2024}. This places machine learning inside a more principled algebraic-geometric framework rather than using it only as a black-box surrogate. The present study shares that motivation: it asks how much of the observed neural improvement is due to geometric structure, as opposed to raw coordinate expressivity.

\subsection{Scientific-ML Context}

The same pattern appears outside Calabi--Yau geometry. The foundational PINN framework treats differential-equation residuals as differentiable training losses~\cite{raissi2019}, and subsequent software and review work has broadened that framework to many forward and inverse PDE problems~\pubcite{6, 37, 38}. At the same time, the PINN optimisation literature makes clear that residual supervision is not automatically enough. Gradient-flow imbalance, neural-tangent-kernel effects, curriculum sensitivity, spectral bias, and activation design can all determine whether a nominally expressive network trains effectively~\pubcite{8--10, 39--41}. Neural-operator and geometric-operator work similarly shows that performance can depend on respecting the structure of the domain and the operator being approximated~\pubcite{7, 11}. Halverson and Ruehle's metric-flow analysis makes a related point in the Calabi--Yau setting: fixed-kernel regimes can be poor metric learners, while finite-width feature learning is important for successful neural metric approximation~\cite{Halverson2024}.

In the specific context of K\"ahler--Einstein metrics, these optimisation pathologies intersect with a deeper issue. The Monge--Amp\`ere equation, patch compatibility, projective consistency, closedness, positivity, and K\"ahler-class conditions are hard geometric constraints, yet they are often enforced through soft penalty terms with hand-tuned weights. Controlling a volume-form residual such as $\omega^3/\Omega\wedge\bar{\Omega}$ also does not by itself prove that the learned two-form lies in the intended K\"ahler class. From an optimisation viewpoint, replacing these constraints by a weighted composite loss relaxes a constrained geometric problem into an unconstrained one, which necessarily permits infeasible solutions during training. The question is therefore not only whether a network can minimise a residual, but whether the ansatz keeps the optimisation close to the geometric constraint manifold.

This broader context shapes the design of the present study. We do not ask whether neural networks can solve the quintic problem in principle; the literature already shows that stronger pipelines can do much better than the simple coordinate baselines used here. We instead ask which part of the representation ladder is responsible for the improvement observed in one archived workflow. The answer is useful precisely because the comparison is local and controlled.

\subsection{Positioning of This Study}

The present paper is therefore closest to an ablation study rather than a new solver paper. It deliberately compares four simplified architecture families under one fixed-$\psi$ protocol, with the specific aim of isolating which geometric properties are responsible for the observed error scale. The SIREN and Cholesky baselines are intentionally simple reference points; they are not claimed to be state-of-the-art representatives of all coordinate or direct-metric methods. Likewise, the regenerated \texttt{CYJAX} Donaldson numbers are same-code algebraic anchors, not like-for-like training-budget competitors.

Within this restricted setting, the comparison isolates one claim: replacing raw coordinate representations by degree-three projective sections changes the attainable accuracy and the stability profile, and this change is the dominant effect in the protocol studied here. The discussion treats atlas methods, spectral models~\pubcite{32, 35}, symmetry canonicalisation, Grassmannian section learning, symbolic constraints, and higher-degree algebraic approximations as motivated next steps rather than defeated alternatives.

\section{Methods}
\label{sec:methodology}

This section defines the geometric setting, the four model architectures, the training objectives, and the experimental design. All quantities reported in Section~\ref{sec:results} are derived directly from the archived configurations and logged checkpoints described here.

\subsection{Geometric Setting and Notation}

We work with the standard Dwork family of quintic hypersurfaces in $\mathbb{CP}^4$~\cite{candelas1991}, defined by the vanishing locus
\begin{equation}
    Q_{\psi}(z)=\sum_{i=0}^{4} z_i^5 - 5 \psi \prod_{i=0}^{4} z_i = 0,
\end{equation}
where $z \in \mathbb{C}^5 \setminus \{0\}$ are homogeneous coordinates, identified modulo overall complex rescaling. The fixed-$\psi$ ablation studied here uses the Fermat point $\psi=0$ throughout as a protocol choice; the extension to $\psi \neq 0$ is left for future work.

For metric evaluation, we select the affine chart determined by the largest homogeneous coordinate in modulus. After permuting indices so that the dominant entry appears first, we introduce affine coordinates
\begin{equation}
    w_i = z_i / z_0, \qquad i=1,\dots,4,
\end{equation}
and write the affine defining equation as
\begin{equation}
    Q_{\psi}^{\mathrm{aff}}(w)=1+\sum_{i=1}^{4} w_i^5 - 5 \psi \prod_{i=1}^{4} w_i = 0.
\end{equation}

\textbf{Notation.} Throughout the paper, $K$ denotes the total K\"ahler potential; $g_{i\bar{j}} = \partial_i \partial_{\bar{j}} K$ is the associated Hermitian metric; $\Omega$ is the holomorphic $(3,0)$-form on the Calabi--Yau threefold; $\eta$ denotes the scalar metric-density used by the potential-based models (defined in Section~\ref{sec:metric-density}); $\kappa$ is its weighted normalisation constant; and $\sigma$ is the weighted root-mean-square deviation of the normalised density from constancy, the primary accuracy diagnostic throughout. The complex-structure parameter $\psi$ is held fixed at zero in all reported experiments.

\subsection{Metric-Density Computation}
\label{sec:metric-density}

The Monge--Amp\`ere equation $\omega^3 = c\,\Omega\wedge\bar\Omega$ is equivalent, after restriction to the hypersurface, to requiring that a certain scalar density---the ratio of the induced metric determinant to the holomorphic volume form---be constant across the manifold. The potential-based models are trained to enforce this constancy. Given a real-valued K\"ahler potential $K$, the ambient Hermitian Hessian is evaluated in affine coordinates by converting the real Hessian to complex form,
\begin{equation}
    H_{i\bar{j}}
    = \frac{1}{4}\left[
        \frac{\partial^2 K}{\partial x_i \partial x_j}
        + \frac{\partial^2 K}{\partial y_i \partial y_j}
        + i\left(
            \frac{\partial^2 K}{\partial x_i \partial y_j}
            - \frac{\partial^2 K}{\partial y_i \partial x_j}
        \right)
    \right],
\end{equation}
where $w_i = x_i + i y_i$. The factor $1/4$ follows from the Wirtinger identities $\partial_{z_i}=\tfrac12(\partial_{x_i}-i\partial_{y_i})$ and $\partial_{\bar z_j}=\tfrac12(\partial_{x_j}+i\partial_{y_j})$, so $H_{i\bar j}$ matches the standard complex Hessian $\partial^2 K/\partial z_i\partial\bar z_j$.

We then form the affine gradient
\begin{equation}
    \frac{\partial Q_{\psi}^{\mathrm{aff}}}{\partial w_i}
    = 5 w_i^4 - 5 \psi \prod_{j \neq i} w_j,
\end{equation}
choose a dependent affine coordinate $\delta$ by the largest gradient component, and pull back the ambient metric to a $3 \times 3$ induced metric on the tangent space. The scalar monitored by the potential-based training path is
\begin{equation}
    \eta(z)
    = \det g_{\mathrm{ind}}(z)\,
      \left|
        \frac{\partial Q_{\psi}^{\mathrm{aff}}}{\partial w_{\delta}}
      \right|^2.
\end{equation}
The factor $|\partial Q_{\psi}^{\mathrm{aff}}/\partial w_\delta|^2$ is the Jacobian of the hypersurface embedding in the chosen affine chart; together with $\det g_{\mathrm{ind}}$, it gives $\eta \propto \det g_{\mathrm{ind}} / |\Omega_{\mathrm{aff}}|^2$, which is precisely the local volume-ratio whose constancy is required by the Monge--Amp\`ere equation. This is the scalar used in both the training loss and the validation $\sigma$ for the SIREN, section-Pure, and section-Hybrid models.

\subsection{Model Architectures}

The fixed-$\psi$ ablation compares four model architectures.

\paragraph{SIREN potential model.}
The coordinate baseline is a SIREN-style periodic-activation network following Sitzmann et al.~\cite{sitzmann2020}. Our implementation uses four layers, width $64$, skip connections, a swept frequency scale $\omega_0$, and a scalar output used as a potential ansatz. It also includes an explicit $\log \|z\|^2$ contribution in the output potential. This Hessian-based construction gives a potential-derived $(1,1)$ form, but it does not by itself guarantee pointwise positivity of the resulting metric throughout training.

\paragraph{Cholesky metric model.}
The positive-definite baseline uses a four-layer SIREN trunk that outputs the parameters of a lower-triangular complex matrix $L$, and the metric is constructed as
\begin{equation}
    g_{\theta}(z)=L_{\theta}(z)L_{\theta}(z)^{\dagger} + 10^{-6} I_3.
\end{equation}
This model bypasses the K\"ahler-potential Hessian and predicts a $3 \times 3$ Hermitian metric directly.

\paragraph{Section-Pure ansatz.}
The algebraic model uses all degree-three monomials in five homogeneous variables, giving $35$ sections,
\begin{equation}
    s_{\alpha}(z)=\prod_{i=0}^{4} z_i^{n_{\alpha,i}},
    \qquad
    \sum_i n_{\alpha,i}=3.
\end{equation}
The learnable object is a Hermitian matrix $H=A A^{\dagger}$, and the model potential is
\begin{equation}
    K_{\mathrm{sec}}(z)=\log \bigl( s(z)^{\dagger} H s(z) \bigr).
\end{equation}
There is no coordinate-dependent neural correction in this variant.

\paragraph{Section-Hybrid ansatz.}
The hybrid model adds a small SIREN correction to the algebraic potential,
\begin{equation}
    K_{\mathrm{hyb}}(z)=K_{\mathrm{sec}}(z) + 0.1\, f_{\theta}(z),
\end{equation}
where $f_{\theta}$ is a four-layer, width-$64$ SIREN network. This is conceptually close to the $\phi$-model family in \texttt{cymetric}~\cite{larfors2022}, with the present-paper modification that the fixed background is the learned section potential rather than the Fubini--Study potential. In this family, the frequency scale $\omega_0$ changes the model class. By contrast, in the pure section model the ablation runner still sweeps $\omega_0$, but that parameter does not enter the algebraic ansatz itself.

\subsection{Section Ansatz and Projective Consistency}

The section-based models require a small amount of explicit geometry to justify the form used in the present study.
The construction is the finite-degree Bergman/Fubini--Study ansatz underlying Donaldson's balanced-metric algorithm and its numerical quintic implementations~\pubcite{13--16, 19, 20}.

\paragraph{Dimension of the degree-three section space.}
The degree-$k$ homogeneous monomials in five variables are counted by stars and bars:
\begin{equation}
    N_k = \binom{k+4}{4}.
\end{equation}
For the archived ablation we use $k=3$, hence
\begin{equation}
    N_3 = \binom{7}{4} = 35,
\end{equation}
which matches the implemented basis size. This is the same low-degree section-ring setting used in numerical quintic metric computations. Douglas et al.~\cite{douglas2008} describe the corresponding quotient construction explicitly: for higher degrees on the hypersurface, the ambient monomials must additionally be quotiented by the defining quintic relation. For example, the regenerated \texttt{CYJAX} $k=6$ Donaldson reference uses $205$ independent hypersurface sections rather than the $210$ ambient degree-six monomials.

\paragraph{Section-Pure potential.}
Let $s(z) \in \mathbb{C}^{N_k}$ denote the vector of degree-$k$ monomials and let
\begin{equation}
    H = A A^{\dagger}
\end{equation}
with $A \in \mathbb{C}^{N_k \times N_k}$ learnable. This guarantees that $H$ is Hermitian positive semidefinite. The model potential is
\begin{equation}
    K_{\mathrm{sec}}(z)=\log \bigl(s(z)^{\dagger} H s(z)\bigr).
\end{equation}
Under the projective rescaling $z \mapsto \lambda z$, every degree-$k$ section transforms as
\begin{equation}
    s(\lambda z)=\lambda^k s(z).
\end{equation}
Therefore
\begin{equation}
    s(\lambda z)^{\dagger} H s(\lambda z)
    =
    |\lambda|^{2k} s(z)^{\dagger} H s(z),
\end{equation}
and hence
\begin{equation}
    K_{\mathrm{sec}}(\lambda z)
    =
    K_{\mathrm{sec}}(z) + 2k \log |\lambda|.
\end{equation}
The additive term depends only on the projective gauge choice, so its $\partial \bar{\partial}$ vanishes on the quotient. Consequently the induced metric derived from $K_{\mathrm{sec}}$ is well defined on projective space. This is the essential mathematical reason the pure section ansatz is geometrically natural.

\paragraph{Section-Hybrid correction.}
For the hybrid model, the scalar potential takes the form
\begin{equation}
    K_{\mathrm{hyb}}(z)=K_{\mathrm{sec}}(z)+\varepsilon f_{\theta}(z),
    \qquad \varepsilon = 0.1.
\end{equation}
At the formal level this preserves the scalar-potential structure,
\begin{equation}
    g^{\mathrm{hyb}}_{i\bar{j}}
    =
    \partial_i \partial_{\bar{j}} K_{\mathrm{hyb}}
    =
    g^{\mathrm{sec}}_{i\bar{j}}
    + \varepsilon\, \partial_i \partial_{\bar{j}} f_{\theta},
\end{equation}
so the hybrid model remains K\"ahler as an ansatz defined by a real potential. What is \emph{not} automatic is projective invariance of the neural correction itself. In the raw implementation, the hybrid branch consumes the flattened real and imaginary parts of the homogeneous representative supplied by the data pipeline. In the projective control variant, it instead consumes the real and imaginary parts of the rank-one projector $z z^\dagger/\|z\|^2$, which is invariant under complex rescaling. We therefore interpret the raw hybrid model conservatively: it is an empirical gauge-fixed correction to the algebraic potential, not a fully derived projective ansatz in the same clean sense as the section-Pure model.

\subsection{Training Objective and Evaluation Metric}

For the potential-based models, the training objective is the weighted variance of $\log \eta$ together with a volume regulariser. Writing
\begin{equation}
    \log \kappa
    = \frac{\sum_i \omega_i \log \eta_i}{\sum_i \omega_i},
\end{equation}
the implemented loss is
\begin{equation}
    \mathcal{L}_{\mathrm{pot}}
    =
    \frac{\sum_i \omega_i (\log \eta_i - \log \kappa)^2}{\sum_i \omega_i}
    + 0.1\,(\log \kappa)^2.
\end{equation}

The Cholesky model uses an analogous loss, but with
\begin{equation}
    \rho_i = \det g_i / |\Omega_i|^2,
\end{equation}
where $|\Omega|^2$ denotes the holomorphic volume-density factor used in the direct-metric baseline. Thus the fixed-$\psi$ ablation does not use a literally identical scalar across all four architectures: three are trained through the induced-metric density $\eta$, whereas the Cholesky baseline uses a direct-metric ratio.

For all four architectures, validation is reported through the RMS dispersion of the normalised monitored scalar,
\begin{equation}
    \sigma
    =
    \sqrt{
        \frac{\sum_i \omega_i (\rho_i / \bar{\rho} - 1)^2}{\sum_i \omega_i}
    },
    \qquad
    \bar{\rho}
    =
    \frac{\sum_i \omega_i \rho_i}{\sum_i \omega_i},
\end{equation}
with $\rho_i = \eta_i$ for the potential-based models and $\rho_i = \det g_i / |\Omega_i|^2$ for the Cholesky baseline. Throughout the paper we report $100 \times \sigma$ as a percentage.

\subsection{Experimental Design}

The fixed-$\psi$ study is a full factorial 96-run ablation. The grid consists of
\begin{itemize}
    \item four model architectures: SIREN, Cholesky, section-Pure, and section-Hybrid;
    \item four learning rates: $\{10^{-5}, 10^{-4}, 10^{-3}, 10^{-2}\}$;
    \item six frequency parameters: $\omega_0 \in \{0.5, 1.0, 1.5, 5.0, 15.0, 30.0\}$.
\end{itemize}
This gives $4 \times 4 \times 6 = 96$ runs.

Each run follows the same architecture and data defaults, together with the architecture-specific overrides required by the ablation design. Concretely, the archived study uses
\begin{itemize}
    \item $\psi = 0$ for all runs;
    \item $100{,}000$ training points and $2{,}000$ validation points;
    \item batch size $20{,}000$ for SIREN and Cholesky;
    \item batch size $2{,}000$ for section-Pure and section-Hybrid;
    \item $200$ training epochs per run;
    \item validation at epoch $1$ and every five epochs thereafter.
\end{itemize}

The data pipeline first generates ambient projective points, projects them to the quintic, and then reprojects the retained points once more before optimisation. The training loop uses geometric weights derived from the projection step and normalises them within each minibatch. In the archived 96-run sweep, the runner does not pass a seed override, so both training and validation data are generated from the same default seed across runs. Consequently, within-archive variation reflects hyperparameter and architecture effects under a fixed dataset rather than sample-to-sample data variation. The supplementary control suite adds explicit seed variation ($11,22,33$).

The JAX training and evaluation scripts used for the archive and the supplementary controls enable 64-bit arithmetic for the geometric computations. We did not, however, run a separate FP32-versus-FP64 ablation, so the precision setting should be treated as part of the fixed protocol rather than as an independently tested factor.

The training loop stores a checkpoint whenever the validation $\sigma$ improves. The reported ``best'' value for a run is therefore the minimum logged validation $\sigma$, not necessarily the final-epoch value. After epoch $50$, training halts early if $\sigma > 100$ or if $\sigma > 2 \sigma_{\mathrm{best}}$ while $\sigma > 5$.

\subsection{Comparability and Scope}

The archive is controlled at the level relevant for an architectural diagnostic: all families target the same fixed quintic, use the same run count, epoch budget, validation cadence, and reporting convention. It is not perfectly symmetric in every implementation detail. The SIREN, section-Pure, and section-Hybrid models are trained through the induced-metric density $\eta$, whereas Cholesky is trained through a direct determinant normalised by the holomorphic volume-density factor. The main archive also uses family-dependent batch sizes and one default seed.

Those asymmetries determine how the results should be read. The 96-run archive is informative about architecture-level behaviour across the stored hyperparameter grid; it is not a full statistical benchmark. The supplementary control suite described in Section~\ref{sec:results-supplementary-controls} addresses the most direct reporting asymmetry by adding explicit seed variation, matched batch sizes, and a shared post-hoc affine-$\eta$ evaluator. The combination supports a local architectural claim, not a universal ranking of Calabi--Yau metric solvers.

\section{Results}
\label{sec:results}

This section summarises the complete fixed-$\psi$ ablation archive used in the present study. For each of the 96 runs, we extract the logged validation $\sigma$ values from the stored training record and report the minimum value attained during that run. A full-length run produces 41 validation measurements: one at epoch $1$ and then one every five epochs through epoch $200$. Because the archived sweep uses a fixed default seed across runs, the results below should be interpreted as a controlled hyperparameter comparison on one shared dataset rather than as a multi-seed statistical estimate. In all results tables, $\sigma$ denotes the weighted RMS relative dispersion of the monitored Monge--Amp\`ere density ratio and is reported as $100\times\sigma$.

\subsection{Ablation Study Overview}

\begin{figure}[t]
\centering
\includegraphics[width=0.88\linewidth]{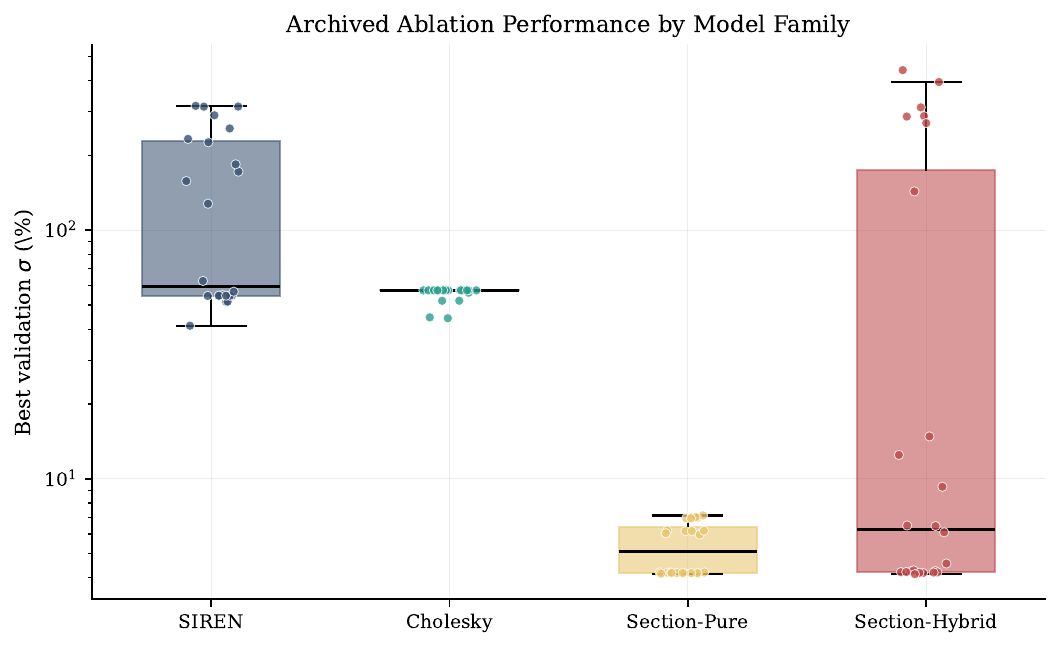}
\caption{Distribution of best validation $\sigma$ (\%) across the 24 archived runs per model architecture. The vertical axis is logarithmic; boxes show medians and interquartile ranges, with whiskers following the standard box-plot convention. Section-based architectures dominate the coordinate-only baselines, while section-Hybrid exhibits a heavier tail of unstable runs than section-Pure.}
\label{fig:family-distributions}
\end{figure}

\begin{figure}[t]
\centering
\includegraphics[width=0.88\linewidth]{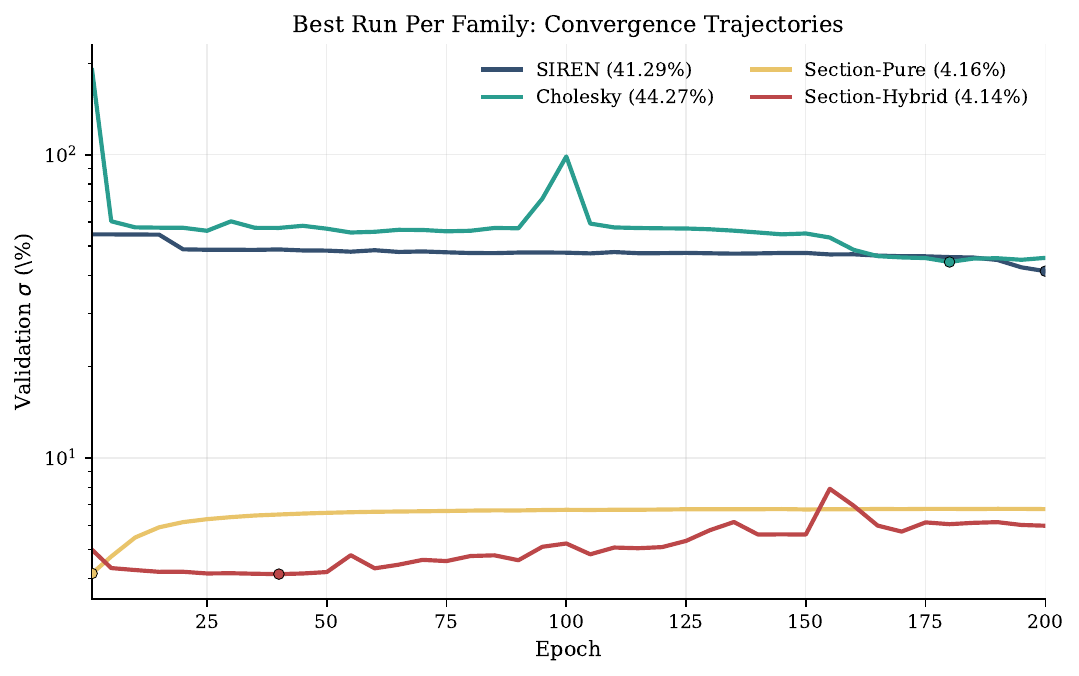}
\caption{Best-run validation $\sigma$ convergence trajectory per model architecture in the fixed-$\psi$ ablation study. Section-based architectures reach the single-digit percentage regime; coordinate-only baselines plateau an order of magnitude higher.}
\label{fig:best-convergence}
\end{figure}

Table~\ref{tab:ablation-summary} reports the central summary statistics for the four model architectures. All values are percentages, i.e.\ $100 \times \sigma$.
Figure~\ref{fig:family-distributions} complements the table with the full distribution of best-of-run outcomes, while Figure~\ref{fig:best-convergence} shows the validation trajectories of the best run per architecture.

\begin{table}[t]
\centering
\caption{Summary statistics for the 96-run fixed-$\psi$ ablation study, grouped by model architecture. Here $\sigma$ is the weighted RMS dispersion of the normalised density ratio. ``Best'' is the minimum validation $\sigma$ attained by any run in the group; ``Median'' and the Q1--Q3 spread are computed over the 24 best-of-run values per architecture.}
\label{tab:ablation-summary}
\begin{tabular}{lccccc}
\hline
Architecture & Runs & Best (\%) & Median (\%) & Q1--Q3 (\%) & Early stop \\
\hline
SIREN & 24 & 41.29 & 59.57 & 54.51--227.78 & 8/24 \\
Cholesky & 24 & 44.27 & 57.31 & 57.29--57.33 & 0/24 \\
Section-Pure & 24 & 4.16 & 5.08 & 4.18--6.38 & 0/24 \\
Section-Hybrid & 24 & 4.14 & 6.28 & 4.22--175.01 & 9/24 \\
\hline
\end{tabular}
\end{table}

We find three main conclusions from the archive.

First, the section-based geometric priors dominate the coordinate-only baselines under this protocol. The best section-Hybrid run achieves $\sigma = 4.14\%$, and the best section-Pure run achieves $\sigma = 4.16\%$, whereas the best SIREN and Cholesky runs remain at $41.29\%$ and $44.27\%$, respectively.

Second, stability differs sharply across families. Section-Hybrid attains the nominally lowest single run in the archive, but it is much less stable than section-Pure: its median best-of-run value is $6.28\%$, its upper quartile reaches $175.01\%$, and 9 of its 24 runs stop early. Section-Pure is effectively tied in best-case performance and is markedly more stable, with no early stops and a median of $5.08\%$.

Third, the two coordinate-based baselines remain far above the section-based models. The direct-metric Cholesky family is stable but limited: its 24 runs occupy a narrow band between $44.27\%$ and $57.40\%$. The SIREN family is less reliable: although its best run reaches $41.29\%$, its median is $59.57\%$, 8 runs stop early, and its worst run reaches $316.57\%$.

\subsection{Best Configurations per Architecture}

\begin{figure*}[t]
\centering
\includegraphics[width=0.96\linewidth]{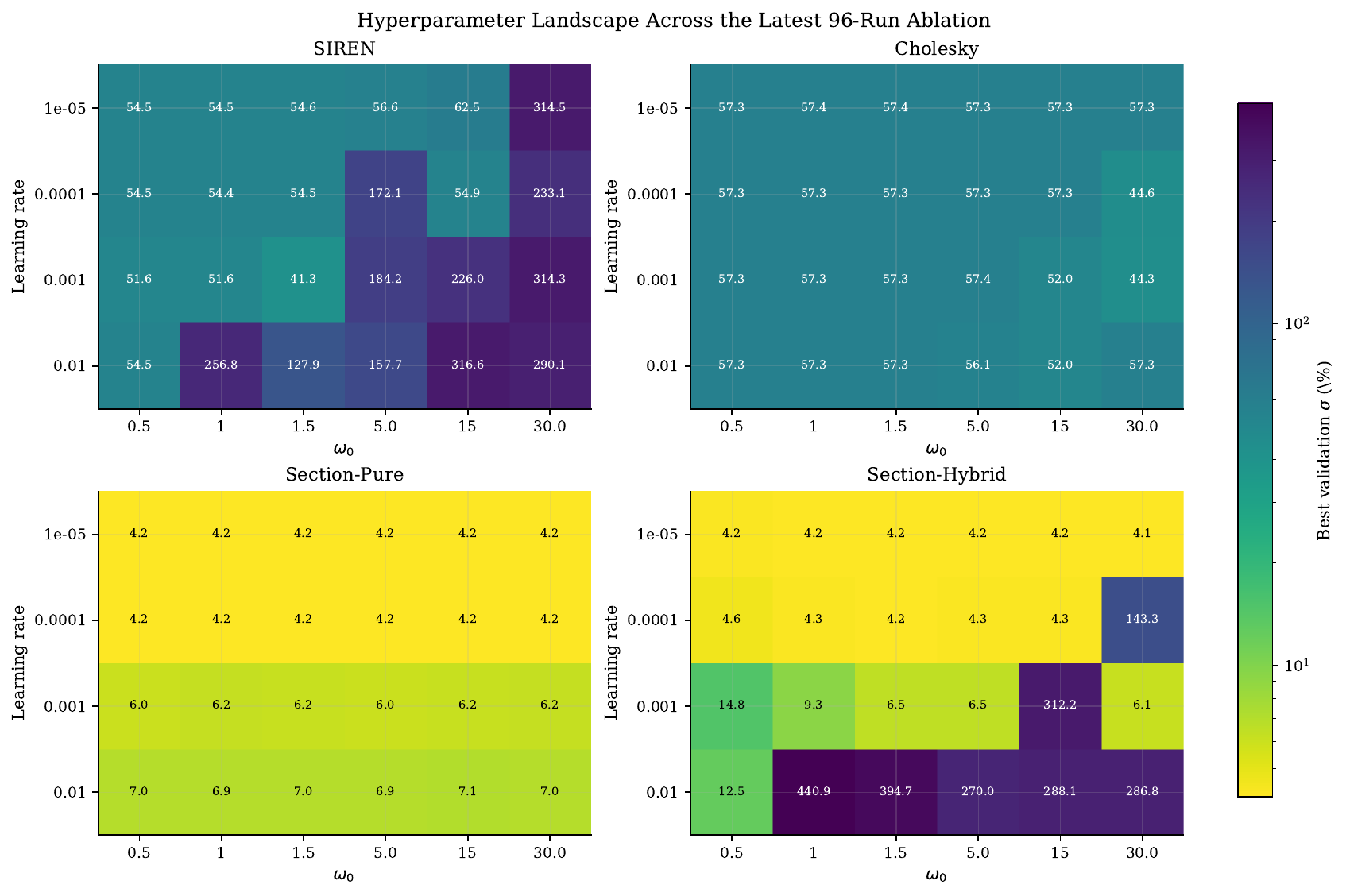}
\caption{Best validation $\sigma$ (\%) for every learning-rate and $\omega_0$ combination in the archived 96-run ablation. The color scale is logarithmic. For section-Pure, the $\omega_0$ sweep is included only because the ablation runner stored it; the pure algebraic model itself does not depend structurally on $\omega_0$.}
\label{fig:ablation-heatmaps}
\end{figure*}

Table~\ref{tab:best-configs} records the best-performing stored run per architecture.
Figure~\ref{fig:ablation-heatmaps} provides the same information at the full grid level by plotting the best attained $\sigma$ for every learning-rate and $\omega_0$ pair.

\begin{table}[t]
\centering
\caption{Best-performing archived run per model architecture. ``Best epoch'' records when the minimum logged validation $\sigma$ was attained; ``Final $\sigma$'' is the value at epoch $200$.}
\label{tab:best-configs}
\begin{tabular}{lccccc}
\hline
Architecture & Learning rate & $\omega_0$ & Best $\sigma$ (\%) & Best epoch & Final $\sigma$ (\%) \\
\hline
SIREN & $10^{-3}$ & 1.5 & 41.29 & 200 & 41.29 \\
Cholesky & $10^{-3}$ & 30.0 & 44.27 & 180 & 45.68 \\
Section-Pure & $10^{-4}$ & 0.5 & 4.16 & 1 & 6.78 \\
Section-Hybrid & $10^{-5}$ & 30.0 & 4.14 & 40 & 5.97 \\
\hline
\end{tabular}
\end{table}

All four family winners reached the full 200-epoch training horizon, so the best configurations were not selected from early-stopped failures. The best checkpoint can nevertheless occur before the final epoch. This matters most for section-Pure: its winning run has $\sigma=4.16\%$ at epoch 1 and $\sigma=6.78\%$ at epoch 200. We therefore do not interpret the section-Pure archive winner as evidence that training monotonically improves the initial algebraic metric. Instead, it shows that the degree-three section representation already begins in the low-error regime, a point reinforced by the fixed identity-$H$ control in Table~\ref{tab:supplementary-controls}.

The section-Pure row requires two cautions. In the current implementation, the pure algebraic model does not depend on $\omega_0$, even though the ablation runner still sweeps that parameter. In addition, the archived sweep uses one fixed default seed and deterministic section initialisation. As a result, the six $\omega_0$ settings at each learning rate are effectively duplicated evaluations, and the family contributes four distinct hyperparameter settings rather than 24 independent stochastic trials.

\subsection{Sensitivity to Hyperparameters}

The heatmaps in Figure~\ref{fig:ablation-heatmaps} show that the four architectures respond to the hyperparameter grid in qualitatively different ways.

For section-Pure, the dominant control parameter is the learning rate rather than $\omega_0$. The two lowest learning rates, $10^{-5}$ and $10^{-4}$, produce a tight $4.16\%$--$4.20\%$ band across all six stored $\omega_0$ values, whereas $10^{-3}$ shifts the family to roughly $6\%$ and $10^{-2}$ to roughly $7\%$. Because $\omega_0$ is inactive in this family, this pattern is best interpreted as a statement about optimisation scale rather than representation scale. In particular, the section ansatz appears robust once the optimiser is placed in the low-learning-rate regime.

Section-Hybrid exhibits a sharper phase transition. At $10^{-5}$, all six runs lie in a narrow and stable $4.14\%$--$4.22\%$ range, making this the only region where the hybrid model is consistently competitive with section-Pure. At $10^{-4}$, five of the six runs remain competitive, but the $\omega_0=30.0$ case already degrades severely to $143.35\%$. At $10^{-3}$ and $10^{-2}$ the family splits between middling runs and outright failures, with several early-stopped trajectories and maxima above $300\%$. This brittle behaviour is consistent with the broader PINN failure-mode literature, where added capacity can worsen optimisation unless the training path is regularised or staged appropriately~\cite{krishnapriyan2021}. The hybrid family therefore shows only a nominal best-case advantage inside a narrow stable band.

The Cholesky and SIREN baselines tell a different story. The Cholesky family is almost frozen near a plateau around $57.3\%$, with only a few visible escape directions: $\omega_0=30.0$ lowers the best runs to the mid-$40\%$ range, and $\omega_0=15$ lowers a smaller set of runs to about $52\%$. The SIREN family is broader and less stable. Low frequencies and moderate learning rates reach the $41\%$--$55\%$ regime, but higher frequencies or aggressive learning rates often trigger degradation into the $150\%$--$316\%$ range. In short, the coordinate-only families remain far less predictable than the section-based ones, and the hybrid family repays its extra flexibility only under careful optimisation.

\subsection{Stability and Failure Analysis}

The archived study also makes it possible to distinguish stable plateaus from brittle best cases. Section-Pure is the clearest example of a stable family: all 24 runs reach epoch 200, none trigger early stopping, and even the worst stored run remains at $7.12\%$. Section-Hybrid has the opposite profile. Although its best run is nominally the lowest in the archive, its worst run reaches $440.86\%$, nine runs stop early, and the upper quartile is already deep in the degraded regime. This is not simply a matter of a heavier tail around a good median; it is evidence that the hybrid branch introduces a real optimisation hazard.

\begin{figure}[t]
\centering
\includegraphics[width=0.88\linewidth]{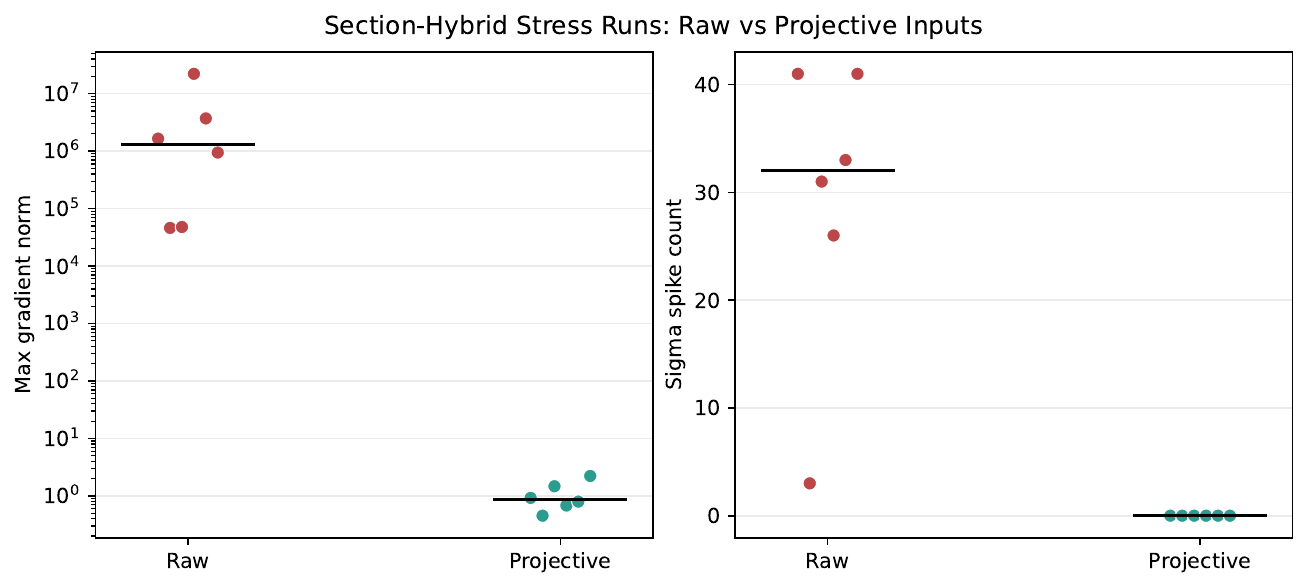}
\caption{Log-derived pathology metrics for the Section-Hybrid stress controls. Each point is one stored stress run. Projective features suppress the extreme gradient and sigma-spike behaviour seen in the raw-coordinate stress runs.}
\label{fig:hybrid-pathology}
\end{figure}

The stress-control logs identify the failure mode more directly. Table~\ref{tab:hybrid-pathology} compares the raw-coordinate and projective-feature hybrid stress runs using only stored training-log quantities: maximum gradient norm, median gradient norm, sigma-spike counts, early stops, and final/best sigma. The raw stress runs show gradient spikes up to $2.2\times 10^{7}$ and accumulate 175 logged sigma spikes above $10\%$, whereas the projective stress runs have no sigma spikes and no early stops. Figure~\ref{fig:hybrid-pathology} visualises the same split. This supports the interpretation that projective features mainly repair optimisation stability rather than producing a clearly lower low-rate optimum.

\begin{table}[t]
\centering
\caption{Optimisation-pathology metrics for the section-Hybrid stress runs. A sigma spike is a logged validation point with $100\times\sigma>10\%$. Gradient norms are read directly from the stored training logs.}
\label{tab:hybrid-pathology}
\begin{tabular}{lcccccc}
\hline
Hybrid input & Runs & Max grad & Median grad & Sigma spikes & Early stop & Final / best $\sigma$ (\%) \\
\hline
Projective features & 6 & 2.2 & 0.087 & 0 & 0 & 6.86 / 4.33 \\
Raw coordinates & 6 & 2.2e+07 & 1.2e+05 & 175 & 2 & 402.72 / 186.03 \\
\hline
\end{tabular}
\end{table}

The distinction between SIREN and Cholesky is subtler. SIREN is both inaccurate and unstable, with eight early stops and a broad spread extending above $300\%$. Cholesky is inaccurate but comparatively controlled: all 24 runs complete, the interquartile band is extremely narrow, and the family occupies a stable plateau rather than a broad failure region. A stable but limited family and an unstable family suggest different next steps, even when neither is competitive with the best section-based models.

\subsection{Control Experiments}
\label{sec:results-supplementary-controls}

To reduce the main cross-architecture asymmetries in the archived comparison, we ran a supplementary control suite with shared affine $\eta = \det g_{\mathrm{loc}} |\partial Q_{\psi}^{\mathrm{aff}} / \partial w_{\delta}|^{2}$, explicit seeds $\{11,22,33\}$, matched batch sizes $\{2000,4096\}$, and a uniform post-hoc evaluator applied to the selected checkpoint of each run. This suite does not replace the 96-run archive; rather, it tests whether the architecture ordering survives under a more standardised evaluation protocol.

\begin{table}[t]
\centering
\caption{Control experiment results under a shared post-hoc affine-$\eta$ evaluator. Here $\sigma$ is the weighted RMS dispersion of the normalised density ratio; all values are percentages and lower is better. Six-run rows correspond to seeds $\{11,22,33\}$ crossed with batch sizes $\{2000,4096\}$.}
\label{tab:supplementary-controls}
\begin{tabular}{lccp{4.8cm}}
\hline
Setting & Runs & $\sigma$ range (\%) & Interpretation \\
\hline
Section-Hybrid (raw, $10^{-5}$) & 6 & 4.02--4.38 & Low-rate raw hybrid; overlaps Section-Pure range \\
Section-Pure & 6 & 4.24--4.39 & Stable algebraic baseline across shared seeds and batch sizes \\
Section-Hybrid (projective, $10^{-5}$) & 6 & 4.25--4.40 & Projective features; no clear best-case gain \\
Section-Hybrid (projective stress, $10^{-4}$) & 6 & 4.26--4.41 & Stress regime remains stable under projective features \\
Section-Hybrid (raw stress, $10^{-4}$) & 6 & 4.45--205.45 & Same stress regime becomes brittle and can fail severely \\
Section identity ($k=3$, $H=\mathbb{I}$) & 1 & 4.44 & Fixed algebraic control without learned refinement \\
Donaldson native ($k=3$) & 1 & 19.38 & \texttt{CYJAX} T-operator balanced-metric reference; not fitted by the present loss \\
Cholesky & 6 & 29.36--29.52 & Stable but limited direct-metric baseline \\
SIREN & 6 & 35.85--36.33 & Weakest learned family under the shared evaluator \\
\hline
\end{tabular}
\end{table}

Three observations are most relevant. First, the coarse architecture ordering survives the standardised evaluation: every section-based setting remains near $4\%$, Cholesky near $29.4\%$, and SIREN near $36.1\%$. Second, the low-learning-rate section-Pure, raw-hybrid, and projective-hybrid ranges overlap strongly, so the controls should not be read as ranking those section-based variants. The robust hybrid conclusion is about stability under stress: the raw stress control spans $4.45\%$--$205.45\%$, whereas the projective stress control stays within $4.26\%$--$4.41\%$. Third, the fixed identity-$H$ degree-three section control already gives $4.443\%$, close to the trained section-Pure and low-learning-rate hybrid models. This indicates that much of the gain comes from the degree-three section prior itself rather than from a large learned correction.

To test whether this ordering is an artefact of a single validation cloud, we then re-evaluated the stored checkpoints on five fresh validation clouds without retraining. Table~\ref{tab:fresh-cloud-robustness} reports the resulting aggregate. Across these clouds, the worst stable section setting remains at least $24.25$ percentage points below the best coordinate-only setting. The deliberately stressed raw hybrid remains pathological and is therefore separated from the stable section settings.

\begin{table}[t]
\centering
\caption{Fresh-cloud robustness check for stored checkpoints. Each checkpoint is re-evaluated on 5 independent validation clouds with 1000 points per cloud. Values are $100\times\sigma$ under the shared affine-$\eta$ evaluator.}
\label{tab:fresh-cloud-robustness}
\begin{tabular}{lccc}
\hline
Setting & Checkpoints $\times$ clouds & Median (\%) & Range (\%) \\
\hline
Section-Hybrid raw low-rate & 30 & 4.29 & 3.84--4.77 \\
Section-Pure & 30 & 4.33 & 4.00--4.71 \\
Section-Hybrid projective low-rate & 30 & 4.35 & 4.02--4.75 \\
Section-Hybrid projective stress & 30 & 4.37 & 4.01--4.77 \\
Section-Hybrid raw stress & 30 & 187.46 & 4.29--212.51 \\
Cholesky & 30 & 29.73 & 28.46--30.59 \\
SIREN & 30 & 36.50 & 34.73--38.13 \\
\hline
\end{tabular}
\par\smallskip\footnotesize Stable section settings beat the best coordinate setting on 5/5 validation clouds; the smallest margin is 24.25 percentage points. The deliberately stressed raw hybrid is reported separately because it is an optimisation-pathology control.
\end{table}

For external context, the same verification pass gives $19.382\%$ for a native \texttt{CYJAX} Donaldson $k=3$ baseline, while a separate native $k=6$ run with $205$ independent hypersurface sections reaches $7.667\%$. The $k=3$ Donaldson comparison uses the same degree-three section space as section-Pure, but it is not a fitted version of the present ansatz: it is the T-operator balanced-metric construction evaluated through its native \texttt{CYJAX} diagnostic. Moreover, the native Donaldson diagnostic is not identical to the weighted RMS statistic used for the learned models. The comparison is therefore a useful algebraic anchor, not a like-for-like training-budget contest.

Table~\ref{tab:anchors-learned-families} reorganises these references into algebraic anchors and learned families. The purpose is to clarify which rows share the affine-$\eta$ evaluator and which rows are included as native Donaldson reference values.

\begin{table}[t]
\centering
\caption{Algebraic anchors and learned-family controls. Values are $100\times\sigma$; learned-family rows report median and range over the six shared-evaluator control checkpoints. Donaldson rows use the native \texttt{CYJAX} diagnostic and are therefore algebraic anchors rather than exact like-for-like entries.}
\label{tab:anchors-learned-families}
\resizebox{\linewidth}{!}{%
\begin{tabular}{lllcp{3.2cm}}
\hline
Class & Method & Diagnostic & $\sigma$ (\%) & Note \\
\hline
Algebraic anchor & Identity section metric ($k=3$, $H=\mathbb{I}$) & shared affine $\eta$ & 4.44 & untrained section ansatz \\
Algebraic anchor & Donaldson balanced metric ($k=3$) & \texttt{CYJAX} native & 19.38 & not the shared affine evaluator \\
Algebraic anchor & Donaldson balanced metric ($k=6$) & \texttt{CYJAX} native & 7.67 & not the shared affine evaluator \\
Learned family & section-Hybrid & shared affine $\eta$ & 4.26 (4.02--4.38) & raw low-rate representative \\
Learned family & section-Pure & shared affine $\eta$ & 4.32 (4.24--4.39) & degree-3 learned section family \\
Learned family & Cholesky & shared affine $\eta$ & 29.43 (29.36--29.52) & direct metric family \\
Learned family & SIREN & shared affine $\eta$ & 36.22 (35.85--36.33) & coordinate potential family \\
\hline
\end{tabular}%
}
\end{table}

\subsection{Empirical Summary}

The archive and controls tell a consistent story. The full grid establishes the hyperparameter landscape: section-Pure is stable, section-Hybrid is powerful but brittle, Cholesky forms a stable high-error plateau, and SIREN is both weaker and less reliable. The control suite then checks the main coarse ranking under a shared evaluator, explicit training seeds, and fresh validation-cloud resampling. In these views, the degree-three section representation is the decisive difference between the low-error and high-error families.

The interpretation of this evidence is taken up in Section~\ref{sec:discussion}. The key empirical point is narrower: the section-versus-coordinate separation is not a single lucky run, but a pattern visible in the distribution of archived runs and preserved by the targeted controls.

\section{Discussion}
\label{sec:discussion}

\subsection{What the Evidence Supports}

The combined evidence supports a local architectural conclusion. In the fixed-$\psi$, 200-epoch, single-patch protocol studied here, section-based models outperform the coordinate-only baselines by roughly an order of magnitude. The archive shows the full hyperparameter landscape: section-Pure and section-Hybrid enter the $4\%$--$6\%$ $\sigma$ regime, while SIREN and Cholesky remain much higher. The supplementary controls then show that this ordering is not removed by explicit seed variation, matched batch sizes, a shared affine-$\eta$ evaluator, or fresh validation-cloud resampling: stable section-based settings remain near $4\%$, Cholesky near $30\%$, and SIREN near $36\%$.

The most conservative interpretation is not that the trained section models discover a radically new correction. The identity degree-three section control already gives $\sigma=4.443\%$, close to the trained section-Pure and low-learning-rate section-Hybrid controls. The dominant gain therefore comes from the projective section representation itself. Training the Hermitian matrix and adding a neural correction provide smaller refinements on top of that structure, and the present data do not statistically separate the low-learning-rate section variants from one another.
This is also consistent with the classical energy-functional lesson of Headrick and Nassar: at fixed section degree, directly optimising a geometric objective can improve over a balanced-metric fixed point, but it does not change the underlying finite-dimensional section space~\cite{headrick2013}. In the present study, the neural ingredient is the residual loss and optimiser wrapped around the degree-three section basis, not a replacement for the algebraic basis itself.

\subsection{Why the Section Prior Helps}

The section ansatz aligns the model class with the geometry of the problem. Degree-three homogeneous sections transform predictably under projective rescaling, so the potential $K_{\mathrm{sec}}=\log(s^\dagger Hs)$ defines a projectively consistent K\"ahler metric after applying $\partial\bar{\partial}$. This built-in structure removes several burdens that a raw coordinate network must learn from data: projective gauge behaviour, holomorphic feature organisation, and compatibility with the algebraic embedding. In a Hessian-based loss, those constraints are not cosmetic; small coordinate-level defects can be amplified by second derivatives.

The Cholesky model illustrates the difference between enforcing one geometric property and encoding the correct geometry. It guarantees pointwise positive definiteness, but it does not guarantee that the metric is K\"ahler, projectively well defined, or derived from a global potential. Its stable plateau is therefore unsurprising: the optimisation finds a controlled direct-metric family, but that family lacks the structural information that makes the section ansatz efficient.
This stable high-error behaviour is qualitatively consistent with the fixed-kernel or weak-feature-learning regimes discussed by Halverson and Ruehle in the Calabi--Yau metric-flow setting~\cite{Halverson2024}; it should be read as a diagnostic analogy, not as a proof that the Cholesky baseline is literally in the same kernel regime.

The hybrid model is more subtle. Adding $f_\theta$ increases expressivity, but the raw correction is evaluated on a chosen representative of the homogeneous coordinates. It therefore reintroduces gauge dependence into an otherwise projective ansatz. This explains the observed split: at gentle learning rates the hybrid branch gives a nominal best archived result, while in the stressed regime it becomes brittle. The projective-feature controls support this interpretation by stabilising the stressed setting without producing a reliable best-case accuracy gain.

\subsection{What Not to Conclude}

The results should not be read as a negative statement about all coordinate-based, direct-metric, neural-operator, atlas, or spectral methods. Stronger pipelines in the literature include ingredients absent from this archive. Examples include moduli-dependent training~\cite{anderson2021}, patch-compatible K\"ahler-potential implementations~\pubcite{25, 26}, symmetry canonicalisation~\cite{hendi_learning_2024}, curvature-aware and phenomenology-oriented metric pipelines~\pubcite{27, 34, 36}, and operator- or spectral-structure-aware methods~\pubcite{11, 32, 35}. The present SIREN and Cholesky baselines are deliberately simple references, not exhaustive representatives of those broader classes.

The study also does not provide formal uncertainty intervals. The archive uses one default seed and family-dependent batch sizes. The control suite varies seeds and batch sizes, and the fresh-cloud pass varies the validation sample, but both are small grids. Their role is to test whether the large ordering in the archive survives obvious protocol asymmetries, not to estimate a full distribution over training data, initialisation, optimiser settings, and evaluator choices. In particular, the $4.14\%$ versus $4.16\%$ archived section-Hybrid/section-Pure gap should be treated as a tie at the current resolution.

Finally, low Monge--Amp\`ere dispersion is not by itself a phenomenology guarantee. Yukawa couplings, kinetic terms, curvature integrals, and Kaluza--Klein spectra depend on inverse metrics, derivatives, harmonic representatives, and global consistency conditions~\pubcite{27, 33--35}. The $\sim4\%$ section-based residuals reported here are best viewed as evidence that the representation is a useful metric-learning substrate and possible warm start for more controlled algebraic, spectral, or atlas-based refinement, not as sufficient precision for downstream string phenomenology.

\subsection{Implications for Future Models}

The ablation points to a concrete design principle: extra flexibility should be added only after the projective and K\"ahler structure is protected. A natural next model is therefore not a larger raw SIREN correction, but a projectively well-defined hybrid branch, a learned section subspace, or a spectral/atlas construction whose features are invariant under the relevant coordinate redundancies.

The next empirical step is similarly clear. A larger frozen-split multi-seed rerun, beyond the checkpoint-only fresh-cloud check reported here, would attach uncertainty intervals to the ordering found here. A fairer direct-metric comparison would train and evaluate all families under a shared scalar diagnostic. Higher-degree section models and on-protocol energy-functional or balanced-metric baselines would determine how much of the remaining gap is due to section degree, optimisation, or evaluator choice. A one-row \texttt{cymetric} $\phi$-model baseline on the same point cloud would also be valuable, because it would connect this diagnostic archive to the de facto neural baseline in the field~\cite{larfors2022}. Symmetry canonicalisation and symbolic constraints are especially natural literature-backed additions~\pubcite{28, 29}. Patch-overlap residuals, a validated curvature evaluator for direct metrics, and an explicit FP32-versus-FP64 precision check are author-proposed protocol checks targeting likely reviewer concerns about coordinate redundancy, differential stability, and numerical precision.

In short, the present archive does not close the Calabi--Yau metric-learning problem. It identifies a narrower mechanism: in this single-patch fixed-$\psi$ setting, projective section structure changes the error scale, while unconstrained neural flexibility mainly changes the stability profile.

\section{Conclusions}
\label{sec:conclusion}

This paper addresses a restricted question: what can be learned from a fixed-$\psi$ archived ablation of single-patch neural and section-based ans\"{a}tze for the quintic Monge--Amp\`ere problem, and does that ordering survive a targeted control experiment pass? Within that scope, the main finding is consistent across both analyses. Section-based geometric priors substantially outperform the coordinate-only SIREN and Cholesky baselines. The archive characterises the architecture-level hyperparameter landscape, while the control experiments recover the same coarse ordering under a shared affine-$\eta$ evaluation: section-based models remain near $4\%$, Cholesky near $29\%$, and SIREN near $36\%$. The degree-three section prior is already strong enough that the fixed identity-$H$ control reaches $4.443\%$; the trained section-Pure architecture remains the most robust learned baseline, and projective hybridisation contributes mainly through improved stability in the stressed regime rather than a reliable peak-accuracy gain.

Geometric structure is therefore not a secondary modification in this protocol; it is the main variable separating the low-error and high-error families. At the same time, the archive supports only a protocol-specific conclusion. The comparison is fixed-$\psi$, single-patch, and tied to the stored workflow. It diagnoses a performance barrier in that regime, but it does not determine what atlas-based, spectral, or properly projective hybrid architectures can achieve.

The practical conclusion is limited but concrete. Within a controlled 96-run fixed-$\psi$ ablation study, supplemented by targeted control experiments and regenerated \texttt{CYJAX} Donaldson references, section-based geometric priors change the observed error scale, coordinate-only models plateau substantially higher, and the more expressive hybrid correction must be stabilised before it can be regarded as a reliable improvement. These findings provide a clear basis for subsequent work on atlas-based, spectral, symmetry-aware, or symbolic architectures.

\section*{Acknowledgments}
Open access funding and article processing charges were supported by Project~DEAKL and the Technical University of Munich (TUM). The author thanks Dr~Johannes Russer (Technical University of Munich, TUM) for helpful feedback on the manuscript structure and presentation. Computations were performed using NVIDIA GPUs and Google Colab infrastructure. The author also gratefully acknowledges the anonymous reviewers for their constructive feedback, which helped improve the clarity and quality of the manuscript.

\section*{Data availability statement}
The data that support the findings of this study are openly available at the following URL/DOI:. \url{10.5281/zenodo.19791880}~\cite{eng2026zenodo}.

\section*{Conflict of interest}
The author declares no competing financial or non-financial interests in relation to the work described.

\section*{Declaration of AI tools usage}
The author used Kimi, ChatGPT, and Claude for code assistance. All content was critically
reviewed and edited by the author. The author takes
full responsibility for the content

\section*{Ethics statement}
This research involves computational and theoretical work only and does not involve human participants, animal subjects, or human data/tissue. No ethics approval was required.

\appendix
\section{Protocol Details and Data Provenance}
\label{app:protocol}

\subsection{Archive Composition}

The analysis in the main text is based on the complete fixed-$\psi$ ablation archive accompanying the project. That archive contains the full $4 \times 4 \times 6$ factorial grid discussed in the main text, together with one additional smoke-test directory. The smoke test is excluded from all summary statistics, tables, and figures reported in this paper. The resulting analytic dataset therefore contains exactly 96 runs.

For each run, the recorded quantity is the validation $\sigma$ evaluated at epoch $1$ and then every five epochs through epoch $200$, unless early stopping occurs first. The best-of-run statistic used throughout the paper is the minimum of that recorded validation sequence. This choice matches the checkpoint-selection logic of the stored training workflow and should be interpreted as a model-selection statistic rather than as a final-epoch statistic.

\subsection{Relation to Earlier Diagnostic Drafts}

Earlier internal diagnostics for this project used a broader act-by-act plateau framing and quoted a degree-three section-Pure plateau near $34\%$ under an older evaluation path. That number is superseded here. The present archive uses the corrected complex-Hessian normalisation, a fixed-$\psi$ 96-run ablation, and the density convention defined in Section~\ref{sec:metric-density}. Consequently, the current $4\%$--$6\%$ section-Pure and section-Hybrid values should be read as results of the corrected fixed-protocol archive, not as directly comparable continuations of the older plateau statistic. The high-error Cholesky behaviour remains part of the story, but it is interpreted more narrowly here as a stable direct-metric baseline plateau within this protocol rather than as the paper's headline lower-bound claim.

\subsection{Extended Summary Statistics}

Table~\ref{tab:appendix-family-summary} supplements the main-text summary by adding the worst stored run and the number of runs that reached epoch $200$.

\begin{table}[t]
\centering
\caption{Extended summary statistics for the 96-run fixed-$\psi$ ablation study, including minimum, quartiles, maximum, and completion counts per model architecture. Here $\sigma$ is the weighted RMS dispersion of the normalised density ratio; all values are percentages.}
\label{tab:appendix-family-summary}
\begin{tabular}{lccccccc}
\hline
Architecture & Min & Q1 & Median & Q3 & Max & Full epoch & Early stop \\
\hline
SIREN & 41.29 & 54.51 & 59.57 & 227.78 & 316.57 & 16/24 & 8/24 \\
Cholesky & 44.27 & 57.29 & 57.31 & 57.33 & 57.40 & 24/24 & 0/24 \\
Section-Pure & 4.16 & 4.18 & 5.08 & 6.38 & 7.12 & 24/24 & 0/24 \\
Section-Hybrid & 4.14 & 4.22 & 6.28 & 175.01 & 440.86 & 15/24 & 9/24 \\
\hline
\end{tabular}
\end{table}

Two features are especially notable. First, the section-based advantage is visible not only in the minima but across the lower and central quantiles of the archived distribution. Second, the distinction between section-Pure and section-Hybrid is most pronounced when one examines the maxima and completion counts rather than the minima alone. The pure section family remains consistently moderate, whereas the hybrid family ranges from very strong to severe failure modes.

\subsection{Validation Curve Details}

The validation curves used for Figure~\ref{fig:best-convergence} are sparse in time compared to the underlying training steps, because the archived workflow validates only periodically. This means that the best-of-run statistic is a minimum over 41 recorded values for a full run, not a continuous minimum over all optimisation steps. That distinction is acceptable for an archive-level paper, but it should be remembered when comparing the present numbers to studies that validate every epoch or every minibatch.

The same caveat applies to early stopping. A run marked as early-stopped has not necessarily diverged monotonically from the start; rather, it violated the stored stopping rule after epoch 50. The early-stop counts should therefore be read as indicators of practical instability under the archived workflow, not as abstract statements about the impossibility of successful training in those parameter regions.

\subsection{Implementation Asymmetries}

Table~\ref{tab:appendix-asymmetries} records the main family-level asymmetries that matter when interpreting the archive.

\begin{table}[t]
\centering
\caption{Implementation asymmetries across the four model architectures. These do not invalidate the architectural comparison but clarify the sense in which the protocol is controlled.}
\label{tab:appendix-asymmetries}
\begin{tabular}{llll}
\hline
Architecture & Monitored scalar & Batch size & Remarks \\
\hline
SIREN & $\eta$ & $20{,}000$ & Potential-based, single-patch \\
Cholesky & $\det g / |\Omega|^2$ & $20{,}000$ & Direct-metric baseline \\
Section-Pure & $\eta$ & $2{,}000$ & $\omega_0$ inactive in model class \\
Section-Hybrid & $\eta$ & $2{,}000$ & Gauge-fixed neural correction \\
\hline
\end{tabular}
\end{table}

This table does not invalidate the comparison, but it clarifies what kind of comparison the paper is making. The archive is controlled in the sense of a shared geometric target and shared high-level workflow, yet it is not perfectly uniform in its internal optimisation mechanics.

\section{Additional Results}
\label{app:additional}

\subsection{Learning-Rate Regimes}

The stored archive is especially informative about learning-rate dependence.

For Section-Pure, the archive separates into three learning-rate bands. The low-learning-rate band ($10^{-5}$ and $10^{-4}$) yields the strongest results, with all runs clustered close to $4.2\%$. The intermediate band ($10^{-3}$) degrades the family to around $6\%$, and the highest band ($10^{-2}$) degrades it further to around $7\%$. Because the family is otherwise so stable, this pattern suggests that the pure section ansatz is limited mainly by optimisation scale rather than by architectural brittleness.

Section-Hybrid is more fragile. At $10^{-5}$ it behaves like a consistently competitive section model. At $10^{-4}$ it is still mostly strong, but a high-frequency outlier has already appeared. At $10^{-3}$ the family becomes mixed, and at $10^{-2}$ most runs are degraded or unstable. This progression is consistent with a hybrid correction that is useful but insufficiently regularised: when optimisation is gentle, the correction can help; when optimisation is aggressive, the added flexibility amplifies instability rather than resolving it.

The coordinate-only families exhibit a weaker and less helpful dependence. Cholesky remains near its plateau for most of the grid, with isolated improvements at specific frequencies. SIREN improves somewhat in the low-frequency, moderate-learning-rate regime, but lacks a comparably clean stable band.

\subsection{Frequency Dependence}

The meaning of the stored frequency parameter differs by family. For SIREN, Cholesky, and section-Hybrid it changes the neural representation directly; for section-Pure it is inert and enters the archive only because the sweep was defined uniformly. This makes section-Pure a useful internal control. The absence of any systematic $\omega_0$ dependence within its low-learning-rate regime confirms that the variation in those rows is not driven by representation changes. With fixed seed and deterministic initialisation in the archived sweep, the six $\omega_0$ entries at each learning rate act as duplicated evaluations for this family.

For the genuinely neural architectures, high frequencies are at best ambiguous and at worst harmful in the archived workflow. SIREN performs best at $\omega_0=1.5$ but becomes highly unstable at $\omega_0=15$ and $30.0$. Cholesky shows a modest improvement window near $\omega_0=30.0$, but otherwise its behaviour is nearly frequency-independent. Section-Hybrid is the most delicate: low learning rates can tolerate large $\omega_0$, but higher learning rates combine destructively with it.

\subsection{Scalar-Curvature Pilot}

As a cautious geometric diagnostic, we evaluated the existing scalar-curvature utility on a 20-point pilot cloud for the best potential-based checkpoint in each family. This is not used as a ranking metric: the calculation differentiates through the metric-density construction, is much more fragile than the Monge--Amp\`ere residual, and the Cholesky baseline is omitted because it outputs a metric directly rather than a K\"ahler potential compatible with this utility. The pilot was retained only because all potential-based rows returned finite values on the sampled points.

\begin{table}[t]
\centering
\caption{Small scalar-curvature pilot for potential-based best checkpoints. This derivative-heavy diagnostic is included only for rows that passed the finite-value screen.}
\label{tab:scalar-curvature-pilot}
\begin{tabular}{lcccc}
\hline
Family & Points & Median $|R|$ & 90th pct. $|R|$ & Max $|R|$ \\
\hline
Section-Hybrid & 20 & 1.93 & 34.9 & 352 \\
Section-Pure & 20 & 1.25 & 2.15 & 4.25 \\
SIREN & 20 & 5.24 & 12.3 & 27.4 \\
\hline
\end{tabular}
\end{table}

The pilot supports the caution in the main text. Section-Pure has small curvature tails on this sample, while section-Hybrid has a finite but much larger tail despite its low Monge--Amp\`ere residual. This is consistent with the view that low $\sigma$ is necessary but not sufficient for downstream curvature-sensitive applications. A larger curvature study would need a dedicated evaluator, patch-aware consistency checks, and a direct-metric curvature path before it should be treated as a primary result.

\subsection{Directions for Future Work}

If this paper were extended beyond the present conservative scope, the natural next empirical additions would be:
\begin{itemize}
    \item a larger multi-seed frozen-protocol rerun, with bootstrap intervals or paired non-parametric tests, to separate hyperparameter effects from seed/data variability;
    \item an on-protocol balanced-metric or energy-functional comparator, including higher section degrees such as $k=6$ or $k=8$;
    \item a \texttt{cymetric} $\phi$-model baseline evaluated on the same point cloud;
    \item a symmetry-aware or canonicalised hybrid model;
    \item patch-overlap residuals and an FP32-versus-FP64 precision check;
    \item a direct comparison between the current section ansatz and a higher-degree section basis.
\end{itemize}

Those additions would not replace the present archive-level analysis; they would test whether the architectural pattern already visible here survives under a stronger experimental standard.

\section*{ORCID iD}
D Yang Eng (\authorchinesename) 0000-0002-0726-1698

% Transcribed from the published journal PDF, DOI 10.1088/2632-2153/ae7494.

% The journal reference list is used here because it differs from the prepublication BibTeX output.

\end{document}